\documentclass[12pt]{amsart}
\usepackage{amsmath,amssymb,amsthm,upref,graphicx,mathrsfs}

\usepackage{color}
\usepackage[
  colorlinks=true,
  linkcolor=blue,
  citecolor=blue,
  urlcolor=blue]{hyperref}

\numberwithin{equation}{section}

\newtheorem{theorem}{Theorem}[section]
\newtheorem{prop}[theorem]{Proposition}
\newtheorem{lemma}[theorem]{Lemma}
\newtheorem{cor}[theorem]{Corollary}

\theoremstyle{definition}
\newtheorem{definition}[theorem]{Definition}

\newtheorem{remark}[theorem]{Remark}

\newcommand{\be}{\begin{eqnarray*}}
\newcommand{\ee}{\end{eqnarray*}}
\newcommand{\beq}{\begin{equation}}
\newcommand{\eeq}{\end{equation}}

\begin{document}

\title[\bf GENERALIZED H\"{O}LDER'S INEQUALITIES AND THEIR APPLICATIONS]
  {\bf H\"{O}LDER'S INEQUALITIES INVOLVING THE INFINITE PRODUCT AND THEIR APPLICATIONS IN MARTINGALE SPACES}

\authors

\author[W. Chen]{Wei Chen}
\address{Wei Chen \\ School of Mathematical Sciences,
Yangzhou University, 225002 Yangzhou, China}
\email{weichen@yzu.edu.cn}

\author[L. Jia]{Longbin Jia}
\address{Longbin Jia \\ School of Mathematical Sciences,
Yangzhou University, 225002 Yangzhou, China}
\email{jialongbin1@126.com}

\author[Y. Jiao]{Yong Jiao}
\address{Yong Jiao\\ School of Mathematics and Statistics, Central South University, 410075 Changsha, China}
\email{jiaoyong@csu.edu.cn}

\makeatletter
\renewcommand{\@makefntext}[1]{#1}
\makeatother \footnotetext{\noindent
The first author supported by the National Natural
Science Foundation of China (Grant No. 11101353), the Natural Science Foundation of Jiangsu
Education Committee (Grant No. 11KJB110018), the Natural Science Foundation of Jiangsu Province (Grant Nos. BK2012682 and BK20141217) and the School Foundation of Yangzhou University (2015CXJ001).
The third author supported by the National Natural
Science Foundation of China (Grant No. 11471337).}


%
%
\begin{abstract} {We give H\"{o}lder's inequalities for integral and conditional expectation involving the infinite product.
Moreover, a generalized Doob maximal operator is introduced and weighted inequalities for
the operator are established.}
\end{abstract}

\maketitle

%
\section{Introduction}
\subsection{Weighted Inequalities for the Hardy-Littlewood Maximal Operator and the Multisubliear Maximal Operator in $R^n$ }
 Let $R^n$ be the $n\hbox{-dimensional}$ real Euclidean
space and $f$ a real valued measurable function. The classical
Hardy-Littlewood maximal operator $M$ is defined by
$$Mf(x)=\sup\limits_{x\in Q}\frac{1}{|Q|}\int_Q|f(y)|dy,$$
where $Q$ is a non-degenerate cube with its sides parallel to the coordinate
axes and $|Q|$ is the Lebesgue measure of $Q.$

Let $u,~v$ be two weights, i.e., positive measurable functions. As is
well known, for $p\geq1,$ Muckenhoupt \cite{Muckenhoupt B} showed that the
inequality
$$ \lambda^p\int_{\{Mf>\lambda\}}u(x)dx
\leq C \int_{R^n}|f(x)|^pv(x)dx, ~~\lambda>0,~f\in{L^p(v)}
$$
holds if and only if $(u,v)\in A_p,$ i.e., for any cube $Q$ in
$R^n$ with sides parallel to the coordinates
$$
\big(\frac{1}{|Q|}\int_Qu(x)dx\big)
\big(\frac{1}{|Q|}\int_Qv(x)^{-\frac{1}{p-1}}dx\big)^{p-1}<C,~p>1;
$$
$$
\frac{1}{|Q|}\int_Qu(x)dx\leq C \mathop{\hbox{ess inf}}\limits_{Q}v(x)
,~p=1.
$$
Suppose that $u=v$ and $p>1,$ Muckenhoupt \cite{Muckenhoupt B} also proved that
$$ \int_{R^n}\big(Mf(x)\big)^pv(x)dx
\leq C \int_{R^n}|f(x)|^pv(x)dx, ~\forall f\in{L^p(v)}
$$
holds if and only if $v$ satisfies
\begin{equation}\label{Ap}
\big(\frac{1}{|Q|}\int_Qv(x)dx\big)
\big(\frac{1}{|Q|}\int_Qv(x)^{-\frac{1}{p-1}}dx\big)^{p-1}<C,~\forall Q.
\end{equation}
The crucial step is to show that if $v$ satisfies the $A_p$,
then there is an $\varepsilon> 0$ such that $v$
also satisfies the $A_{p-\varepsilon}.$
But, the problem of finding all $u$ and $v$ such that
$$ \int_{R^n}\big(Mf(x)\big)^pu(x)dx
\leq C \int_{R^n}|f(x)|^pv(x)dx,~\forall f\in{L^p(v)}
$$
is much more complicated. In order to solve the problem,
Sawyer \cite{Sawyer E T.} established the
testing condition $S_{p,q},$ i.e., for any cube $Q$ in
$R^n$ with sides parallel to the coordinates
$$\Big(\int_{Q}\big(M(\chi_Qv^{1-p'})(x)\big)^qu(x)dx\Big)^{\frac{1}{q}}
\leq C\big(\int_Qv(x)^{1-p'}dx\big)^{\frac{1}{p}},~\forall Q$$
where $1<p\leq q<\infty.$
The condition $S_{p,q}$ is a
sufficient and necessary condition such that
the weighted inequality
$$\Big(\int_{R^n}\big(Mf(x)\big)^qu(x)dx\Big)^{\frac{1}{q}} \leq C
\Big(\int_{R^n}|f(x)|^pv(x)dx\Big)^{\frac{1}{p}}, ~\forall f\in{L^p(v)}$$
holds.
In this case, the method of proof is very interesting. Motivated by these
results, the theory of weighted inequalities developed rapidly in the last few decades
, not
only for the Hardy-Littlewood maximal operator but also for some of the main
operators in Harmonic Analysis like Calders\'{o}n-Zygmund operators (see \cite{Garcia Rubio}
and \cite{Cruz-Uribe D. J. M. Martell} for more information).

Recently, the multisublinear maximal function
\begin{equation}\label{multi_maximal_operator}\mathcal{M}(f_1,...,f_m)(x) = \sup\limits_{x\in Q}
\prod\limits_{i=1}\limits^{m}\frac{1}{|Q|}\int_Q|f_i(y_i)|dy_i
\end{equation}
associated with cubes with sides parallel to the coordinate
axes was studied in \cite{Lerner A.K. Ombrosi S.}.
The importance of this operator is that it generalizes the Hardy--Littlewood
maximal function (case $m=1$) and in several ways it controls the class
of multilinear Calder\'{o}n--Zygmund operators as it was shown in \cite{Lerner A.K. Ombrosi S.}.
The relevant class of multiple weights for $\mathcal{M}$ is given by the condition $A_{\overrightarrow{p}}:$ for
$\overrightarrow{p}=(p_1,p_2,\cdot\cdot\cdot,~ p_m),$
$\overrightarrow{\omega}=(\omega_1, ~\omega_2,\cdot\cdot\cdot,~\omega_m)$ and a weight $v,$
the weight vector $(v, \overrightarrow{\omega})\in A_{\overrightarrow{p}}$ if
$$\sup_Q\frac{v(Q)}{|Q|}\prod\limits^m_{i=1}\big(\frac{\sigma_i(Q)}{|Q|}\big)^{\frac{p}{p'_i}}
< \infty,$$
where $\frac{1}{p}=\sum\limits^m_{i=1}\frac{1}{p_i }$ and $1\leq p_1,p_2,...,p_m<\infty.$

It is easy to see that in the linear case (that is,
if $m=1$), condition $A_{\overrightarrow{p}}$ is the usual $A_p.$
In \cite{Lerner A.K. Ombrosi S.} the following multilinear extension
of the Muckenhoupt $A_p$ theorem for the maximal function was obtained. The inequality
$$\|\mathcal{M}(\overrightarrow{f})\|_{L^{p,\infty}(v)}\leq
C\prod\limits^m_{i=1}\|f_i\|_{L^{p_i}(\omega_i)},
~\forall f_i\in L^{p_i}(\omega_i)$$
holds if and only if $(v, \overrightarrow{\omega})\in A_{\overrightarrow{p}}.$
Moreover, if $1< p_1,p_2,...,p_m<\infty$ and
$v=\prod_{i=1}^mw_i^{p/p_i},$ then the multilinear $A_{\overrightarrow{p}}$ condition
has the characterization in terms of the linear $A_p$ classes, i.e.,
\begin{equation}\label{cha_Ap}(v, \overrightarrow{\omega})\in A_{\overrightarrow{p}} \hbox{ if and only if } v\in A_{mp}\hbox{ and }\omega_j^{1-p_j'}\in A_{mp_j'},~j=1,...,m.\end{equation}
Employing the characterization, they got that the inequality
$$\|\mathcal{M}(\overrightarrow{f})\|_{L^{p}(v)}\leq
C\prod\limits^m_{i=1}\|f_i\|_{L^{p_i}(\omega_i)},
~\forall f_i\in L^{p_i}(\omega_i)$$
holds if and only if $(v, \overrightarrow{\omega})\in A_{\overrightarrow{p}}.$ The more general case was extensively discussed in
\cite{Grafakos Liu, Grafakos Perez}. Recently, Dami\'{a}n, Lerner and P\'{e}rez \cite{W. Damian} observed that \begin{equation}\label{Con_dy}\mathcal{M}(\overrightarrow{f})\leq 6^{mn}\sum_{\alpha=1}^{2^n}\mathcal{M}^{\mathcal{D}_{\alpha}}(\overrightarrow{f}).\end{equation}
Using the observation, they obtained a sharp mixed $A_p-A_\infty$ bound for the operator.
In order to establish the generalization of Sawyer's theorem to the multilinear setting,
a kind of monotone property and a
reverse H\"{o}lder's inequality on the weights were introduced in \cite{W. M. Li} and \cite{Chen-Damian}, respectively.
They both established the multilinear version of Sawyer's result.
Recently, Li and Sun \cite{Li_Sun} made progress for $\mathcal{M}$ by \eqref{Con_dy}.
Moreover, the multilinear fractional maximal operator and the multilinear fractional
strong maximal operator associated with rectangles were studied in \cite{cao-xue} and \cite{cao-xue-yabuta}, respectively.
The new methods, including atomic decomposition of tent space in \cite{cao-xue} and Carleson embedding theorem in
\cite{cao-xue-yabuta}, may provide different approaches to deal with two-weight norm inequalities.

In this paper, we define a new generalized maximal function
\begin{equation*}\label{g-multi_maximal_operator}\mathfrak{M}(\overrightarrow{f})(x)\triangleq\sup\limits_{x\in Q}
\prod\limits_{i=1}\limits^{\infty}\frac{1}{|Q|}\int_Q|f_i(y_i)|dy_i
\end{equation*}
for suitable $\overrightarrow{f}=(f_1,f_2,...)$(see Lemma \ref{HL-constant} for a kind of suitable
condition).
Then it is natural to establish weighted inequalities for it.
Unfortunately, methods of \cite{Lerner A.K. Ombrosi S.,W. Damian} are not suitable.
One reason is that \eqref{cha_Ap} and \eqref{Con_dy} are invalid when $m=\infty.$
However, this is not the end of the story. The first author \cite{Chen_Zhang}
defined a generalized dyadic maximal operator involving the infinite product
and discussed weighted inequalities for the operator
by a formulation of the Carleson embedding theorem. Now,
we will establish related theory in martingale setting.

\subsection{Weighted Inequalities for the Doob Maximal Operator and the Multisubliear Doob Maximal Operator in Martingale Setting }

Let $(\Omega,\mathcal {F},\mu)$ be a complete probability space and let
$(\mathcal {F}_n)_{n\geq0}$ be an increasing sequence of
sub$\hbox{-}\sigma\hbox{-}$fields of $\mathcal{F}$ with
$\mathcal{F}=\bigvee\limits_{n\geq0}\mathcal{F}_n.$ A weight
$\omega$ is a random variable with $\omega>0$ and $
E(\omega)<\infty.$ For any $n\geq0$ and integral function $f,$ we denote the conditional expectation
with respect to $\mathcal {F}_n$ by $E_n(f)$ or $E(f|\mathcal {F}_n),$ then $(E_n(f))_{n\geq0}$ is
an uniformly integral martingale. For $(\Omega,\mathcal {F},\mu)$ and $(\mathcal
{F}_n)_{n\geq0},$ the family of all stopping times is denoted by
$\mathcal {T}.$ Given $\tau\in \mathcal {T},$ let
$$\mathcal{F}_{\tau}=\{F\in\mathcal{F}:~F\cap\{\tau\leq n\}\in \mathcal{F},~\forall n\geq0 \},$$
then $\mathcal{F_{\tau}}$ is a sub$\hbox{-}\sigma\hbox{-}$field of $\mathcal{F}.$
For an integral function $f,$ we denote the conditional expectation
with respect to $\mathcal {F}_{\tau}$ by $E_{\tau}(f).$ Moreover, if we define
$f_{\tau}(x)\triangleq f_{\tau(x)}(x)\chi_{\{\tau<\infty\}}+f(x)\chi_{\{\tau=\infty\}},$
then $E_{\tau}(f)=f_{\tau}$ (see \cite{Neveu,R. L. Long} for more information).
Let $B\in \mathcal {F},$
we always denote $\int_\Omega\chi_Bd\mu$ and
$\int_\Omega\chi_B\omega d\mu$ by $|B|$ and $|B|_\omega,$
respectively.

Suppose that functions $f,~g$ are
integrable on the probability space $(\Omega,\mathcal {F},\mu),$ then the Doob maximal operator
and the bilinear Doob maximal operator are defined by
$$Mf=\sup_{n\geq 0}|E_n(f)|\text{ and }
\mathcal{M}(f,g)=\sup_{n\geq 0}|E_n(f)||E_n(g)|,$$
respectively.

In regular martingale spaces, Izumisawa and Kazamaki \cite{M.
Izumisawa} characterized the inequality
$$\Big(\int_\Omega(Mf)^pvd\mu\Big)^{\frac{1}{p}} \leq C
\Big(\int_\Omega|f|^pvd\mu\Big)^{\frac{1}{p}},$$ where $p>1$ and $v$
is a weight. In addition, Long and Peng \cite{Long R. L. and Peng L. Z.} obtained
probabilistic $A_p$ condition and $S_p$ condition,
which were also discussed in \cite{M. Kikuchi} and \cite{X. Q. Chang}, respectively.

Let $v,~\omega_1, ~\omega_2$ be weights and $1<p_1, ~p_2<\infty.$
Suppose that
$\frac{1}{p}=\frac{1}{p_1 }+\frac{1}{p_2 }$ and $(\omega_1, \omega_2)\in RH(p_1, p_2),$
for the bilinear Doob maximal operator $\mathcal{M}$, Chen and Liu \cite{W. Chen}
characterized the weights for which $\mathcal{M}$ is bounded from
$L^{p_1}(\omega_1)\times L^{p_2}(\omega_2 )$
to $L^{p,\infty}(v)\hbox{ or }L^p(v).$ If $v=\omega_2^{\frac{p}{p_2}}\omega_2^{\frac{p}{p_2}},$
they also have a bilinear version for the convergence of martingale.

In this paper, we define the generalized Doob maximal operator $\mathfrak{M}$
in the following way:
$$\mathfrak{M}(\overrightarrow{f})\triangleq\sup\limits_{n\geq0}\prod\limits^{\infty}_{i=1}|E_n(f_i)|,$$
where $\overrightarrow{f}=(f_1,f_2,...)$ and $\overrightarrow{f}$ is subjected to suitable restrictions. The suitable
restrictions can be found in Proposition \ref{C-P-lem-def} and Remark \ref{C-P-lem-finite}.
Now, we state our main results.
\begin{theorem} \label{thm_AP}Let $v$ be a weight and $\overrightarrow{\omega}\in RH_{\overrightarrow{p}},$
then the following statements
are equivalent:
\begin{enumerate}
\item \label{C-L-P-W-1}There exists a positive constant $C$ such that
      \begin{equation}\label{Th_B_1}\Big(\int_{\{\tau<\infty\}}
      \prod\limits_{i=1}^{\infty}E_{\tau}(f_i)^pvd\mu\Big)^{\frac{1}{p}}
      \leq C\prod\limits_{i=1}^{\infty}\|f_i\|_{L^{p_i}(\omega_i)},
      ~\forall \tau\in\mathcal{T},~ f_i\in L^{p_i}(\omega_i),~i\in N,\end{equation}
where $\prod\limits_{i=1}^{\infty}\|f_i\|_{L^{p_i}(\omega_i)}<\infty;$
\item \label{C-L-P-W-2} There exists a positive constant $C$ such that
      \begin{equation}\label{Th_B_2}\|\mathfrak{M}(\overrightarrow f)\|_{L^{p,\infty}(v)}\leq
       C\prod\limits_{i=1}^{\infty}\|f_i\|_{L^{p_i}(\omega_i)},
      ~\forall f_i\in L^{p_i}(\omega_i), ~i\in N,\end{equation}
where $\prod\limits_{i=1}^{\infty}\|f_i\|_{L^{p_i}(\omega_i)}<\infty;$
\item \label{C-L-P-W-3}The weight vector $(v,\overrightarrow{\omega})$
satisfies the condition $A_{\overrightarrow{p}},$ i.e.,
\begin{equation}\label{Th_B_3}
(v,\overrightarrow{\omega})\in A_{\overrightarrow{p}}.
      \end{equation}
\end{enumerate}
\end{theorem}

\begin{theorem}\label{thm_Sp} Let $v$ be a weight and $\overrightarrow{\omega}\in RH_{\overrightarrow{p}},$
then the following statements
are equivalent:
\begin{enumerate}
\item \label{C_L_P_thm_Sp_1}There exists a positive constant $C$ such that
\be \|\mathfrak{M}(\overrightarrow f)\|_{L^p(v)}\leq
       C\prod\limits_{i=1}^{\infty}\|f_i\|_{L^{p_i}(\omega_i)},
      ~\forall f_i\in L^{p_i}(\omega_i),~i\in N,
\ee
where $\prod\limits_{i=1}^{\infty}\|f_i\|_{L^{p_i}(\omega_i)}<\infty;$
\item \label{C_L_P_thm_Sp_2}There exists a positive constant $C$ such that
\begin{equation}
\label{Th_A_1}
\|\mathfrak{M}(\overrightarrow {g\sigma})\|_{L^p(v)}\leq
       C\prod\limits_{i=1}^{\infty}\|g_i\|_{L^{p_i}(\sigma_i)},
      ~\forall g_i\in L^{p_i}(\sigma_i),~i\in N,
\end{equation}
where $\prod\limits_{i=1}^{\infty}\|g_i\|_{L^{p_i}(\sigma_i)}<\infty;$
\item \label{C_L_P_thm_Sp_3} The weight vector $(v,\overrightarrow{\omega})$
satisfies the condition $S_{\overrightarrow{p}},$ i.e.,
\be
(v,\overrightarrow{\omega})\in S_{\overrightarrow{p}}.
\ee  \end{enumerate}
\end{theorem}

The remainder of this paper is organized as follows. In Section \ref{section B},
we prove the generalized H\"{o}lder inequalities for integral and conditional expectation in details,
which will be used in Section \ref{section C}.
The proofs of Theorem \ref{thm_AP}
and Theorem \ref{thm_Sp} are contained in Section \ref{section C}.
In this paper, for simplicity, we omit the annotation `almost everywhere' in the following statements.

\section{Generalized H\"{o}lder Inequalities for Integral and Conditional Expectation}\label{section B}

The section consists of a series of Lemmas. If the readers are familiar with them,
they may read ahead to Theorems \ref{C-P-prop-A}, \ref{C-P-series-cor-E} and \ref{C-P-cor-gen-h} directly.

\subsection{Some Properties of Series, Lebesgue's Integral and Infinite Product}\label{section B0}
Let $\{a_i\}$ be a sequence of real numbers.
Let $\{s_n\}$ be the sequence obtained from $\{a_i\},$ where for each $n\in N,~s_n=\sum\limits^{n}_{i=1}a_i.$
If $s_n$ converges in $R$ or diverges to $+\infty$ (or $-\infty$),
we say that the sum of the series is well defined and we denote the sum as $\sum\limits^{\infty}_{i=1}a_i.$
Let $\lambda_i\in(0,1),b_i\in R, ~i\in N,$ and let $\sum\limits_{i=1}^{\infty}\lambda_i=1.$ It is known that $(N, 2^N)$ is a measurable space.
By the sequences $\{\lambda_i\}$ and $\{b_i\},$
we can define a measure $\lambda$ and a measurable function $b$ on the space in the following way
$$\lambda(i)=\lambda_i\hbox{ and }b(i)=b_i,~\forall i\in N.$$ Then $(N, 2^N, \lambda)$ is a probability space. Applying Levi's Lemma, we have
$$\sum\limits^{\infty}_{i=1}\lambda_ib^+_i=\lim\limits_{k\rightarrow\infty}\sum\limits^k_{i=1}\lambda_ib^+_i
=\lim\limits_{k\rightarrow\infty}\int_N b^+\chi_{\{1,2,...,k\}} d\lambda=\int_N b^+d\lambda$$ and $$\sum\limits^{\infty}_{i=1}\lambda_ib^-_i=\lim\limits_{k\rightarrow\infty}\sum\limits^k_{i=1}\lambda_ib^-_i
=\lim\limits_{k\rightarrow\infty}\int_N b^-\chi_{\{1,2,...,k\}} d\lambda=\int_N b^-d\lambda.$$
For simplicity, we denote $\sum\limits^{\infty}_{i=1}\lambda_ib^+_i$ and $\sum\limits^{\infty}_{i=1}\lambda_ib^-_i$ by $A$ and $B,$ respectively.  It follows that $A,B\in [0,+\infty].$ If $A$ or $B$ is finite,
then $\sum\limits^{\infty}_{i=1}\lambda_ib_i$ is well defined, integral of $b$ exists and
$$\sum\limits^{\infty}_{i=1}\lambda_ib_i=\int_N bd\lambda.$$

This paper also involves the concept of an infinite product. Let us recall the definition (see, e.g., \cite[p.~298]{Rudin}).

\begin{definition} Suppose $\{c_n\}$ is a sequence of complex number,
$$p_n=\prod\limits^{n}_{i=1}c_i,$$
and $p=\lim\limits_{n\rightarrow\infty}p_n$ exists. Then we write
\begin{equation}\label{Rudin}p=\prod\limits^{\infty}_{i=1}c_i.\end{equation}
The $p_n$ are the partial products of the infinite product (\ref{Rudin}).
We should say that the infinite product (\ref{Rudin}) converges if the sequence $\{p_n\}$ converges.
\end{definition}

\begin{remark} \label{re-aa}Suppose $\{c_n\}$ and $\{c'_n\}$ are nonnegative sequences,
and the infinite product $\prod\limits^{\infty}_{i=1}c_i$ converges. If $c'_n\leq c_n,~n\in N,$ then the infinite product $\prod\limits^{\infty}_{i=1}c'_i$ also converges.
\end{remark}

\begin{remark} \label{re-bb}Suppose $\{f_i\}$ is a sequence of measurable functions on a measurable space $(\Omega,\mathcal{F}),$ and suppose that
the sequence of numbers $\{\prod\limits^{n}_{i=1}f_i(x)\}$ converges for every $x\in \Omega.$ We can then define a function $\prod\limits^{\infty}_{i=1}f_i$ by
$$\prod\limits^{\infty}_{i=1}f_i(x)=\lim\limits_{n\rightarrow\infty}\prod\limits^{n}_{n=1}f_i(x).$$
Thus the function $\prod\limits^{\infty}_{i=1}f_i(x)$ is well defined.
\end{remark}

\begin{lemma}\label{C-P-cor-Jensen}
Let $(\Omega,\mathcal{F},\mu)$ be a probability space. If the measurable function $f:\Omega\rightarrow R$ such that $\exp(f)$ is integrable, then integral of the function f exists and $$\exp\big(\int_\Omega f d\mu\big)\leq\int_\Omega\exp(f)d\mu.$$
\end{lemma}
\noindent{\bf Proof of Lemma \ref{C-P-cor-Jensen}} It is clear that $f^+\leq \exp(f^+)=\max\{\exp(f),1\}\leq\exp(f)+1,$ then
$$\int_\Omega f^+d\mu\leq\int_\Omega \exp(f)d\mu+1<\infty.$$
Thus integral of the measurable function $f$ exists. If $\int_\Omega f^-d\mu<\infty,$ it follows from Jensen's inequality that $$\exp\big(\int_\Omega f d\mu\big)\leq\int_\Omega\exp(f)d\mu.$$ If $\int_\Omega f^-d\mu=+\infty,$ we have $\int_\Omega fd\mu=-\infty$ and $\exp\big(\int_\Omega f d\mu\big)=0.$ We are done.$\blacksquare$

\begin{cor}\label{C-P-lem-B}
Let $\lambda_i\in(0,1), i\in N,$ $\sum\limits_{i=1}^{\infty}\lambda_i=1.$ If $b_i\in R,~i\in N$ and $\sum\limits^{\infty}_{i=1}\lambda_i\exp (b_i)<\infty,$ then $\sum\limits^{\infty}_{i=1}\lambda_ib_i$ is well defined and
$$\exp(\sum\limits^{\infty}_{i=1}\lambda_ib_i)\leq \sum\limits^{\infty}_{i=1}\lambda_i\exp (b_i).$$
\end{cor}
\noindent{\bf Proof of Corollary \ref{C-P-lem-B}} The corollary is another version of Lemma \ref{C-P-cor-Jensen}. We can prove the corollary in the way of Lemma \ref{C-P-cor-Jensen} with obvious changes and we omit it.$\blacksquare$

\begin{lemma}\label{C-P-lem-C}
Let $\lambda_i\in(0,1), i\in N$ and $\sum\limits_{i=1}^{\infty}\lambda_i=1.$ If $a_i\geq0,~i\in N$ and $\sum\limits^{\infty}_{i=1}\lambda_ia_i<\infty,$  then
$$\prod\limits^{\infty}_{i=1}a_i^{\lambda_i}\leq \sum\limits^{\infty}_{i=1}\lambda_ia_i.$$
\end{lemma}

\noindent{\bf Proof of Lemma \ref{C-P-lem-C}} Without loss of generalization, we assume $a_i>0, ~i\in N.$
Substituting $b_i=\ln a_i, ~ i \in N$ into Corollary \ref{C-P-lem-B}, we have
$$\exp(\sum\limits^{\infty}_{i=1}\lambda_i\ln a_i)\leq \sum\limits^{\infty}_{i=1}\lambda_i\exp (\ln a_i).$$
It follows that
$$\prod\limits^{\infty}_{i=1}a_i^{\lambda_i}\leq \sum\limits^{\infty}_{i=1}\lambda_ia_i.\blacksquare$$

\begin{lemma}\label{C-P-lem-D}
Let $1<p_i<\infty, i\in N$ and $\sum\limits_{i=1}^{\infty}\frac{1}{p_i}=1.$ If $a_i\geq0,~i\in N$ and $\sum\limits^{\infty}_{i=1}\frac{a_i}{p_i}<\infty,$  then
$$\prod\limits^{\infty}_{i=1}a_i^{\frac{1}{p_i}}\leq \sum\limits^{\infty}_{i=1}\frac{a_i}{p_i}.$$
\end{lemma}

\noindent{\bf Proof of Lemma \ref{C-P-lem-D}} Substituting $\lambda_i=\frac{1}{p_i},
~i \in N$ into Lemma \ref{C-P-lem-C}, we have Lemma \ref{C-P-lem-D}.$\blacksquare$

\begin{lemma}\label{C-P-lem-E}
Let $1<p_i<\infty, i\in N$ and $\sum\limits_{i=1}^{\infty}\frac{1}{p_i}=1.$  If $c_i\geq0,~i\in N$ and $\sum\limits^{\infty}_{i=1}\frac{c_i^{p_i}}{p_i}<\infty,$ then
$$\prod\limits^{\infty}_{i=1}c_i\leq \sum\limits^{\infty}_{i=1}\frac{c_i^{p_i}}{p_i}.$$
\end{lemma}

\noindent{\bf Proof of Lemma \ref{C-P-lem-E}} Substituting $a_i=c_i^{p_i}, ~i \in N$
into Lemma \ref{C-P-lem-D}, we have Lemma \ref{C-P-lem-E}.$\blacksquare$

\subsection{Generalized H\"{o}lder's Inequality for Integral}\label{section B 1}
In the subsection,
we suppose that $(\Omega,\mathcal{F},\mu)$ is a measure space and $\{f_i\}$ is a sequence of nonnegative measurable functions on $(\Omega,\mathcal{F},\mu)$.

\begin{lemma}\label{C-P-lem-F}
Let $1<p_i<\infty\hbox{ and }\|f_i\|_{L^{p_i}}=1,~ i\in N.$ If $\sum\limits_{i=1}^{\infty}\frac{1}{p_i}=1,$
then the function $\prod\limits^{\infty}_{i=1}f_i$ is well defined and $$\|\prod\limits^{\infty}_{i=1}f_i\|_{L^1}\leq 1.$$
\end{lemma}

\noindent{\bf Proof of Lemma \ref{C-P-lem-F}} Since $\|f_i\|_{L^{p_i}}=1,~ i\in N$
and $\sum\limits_{i=1}^{\infty}\frac{1}{p_i}=1,$ we have
$$\int_{\Omega}\sum\limits_{i=1}^{\infty}\frac{{f_i}^{p_i}}{p_i}d\mu
=\sum\limits_{i=1}^{\infty}\int_{\Omega}\frac{{f_i}^{p_i}}{p_i}d\mu
=\sum\limits_{i=1}^{\infty}\frac{\int_{\Omega}{f_i}^{p_i}d\mu}{p_i}
=\sum\limits_{i=1}^{\infty}\frac{1}{p_i}=1<\infty,$$
where we have used the monotone convergence theorem.
It follows that $$\sum\limits_{i=1}^{\infty}\frac{{f_i}^{p_i}}{p_i}<\infty.$$
Combining this with Lemma \ref{C-P-lem-E}, we get that $\prod\limits^{\infty}_{i=1}f_i$ is well defined and
$$\prod\limits^{\infty}_{i=1}f_i\leq \sum\limits^{\infty}_{i=1}\frac{f_i^{p_i}}{p_i}<\infty.$$
Hence, $$\int_{\Omega}\prod\limits^{\infty}_{i=1}f_id\mu\leq \sum\limits^{\infty}_{i=1}\frac{\int_{\Omega}f_i^{p_i}d\mu}{p_i}
=\sum\limits^{\infty}_{i=1}\frac{1}{p_i}=1.\blacksquare$$

\begin{lemma}\label{C-P-lem-G}
Let $1<p_i<\infty,~ i\in N$ and $\sum\limits_{i=1}^{\infty}\frac{1}{p_i}=1.$ If $\prod\limits^{\infty}_{i=1}\|f_i\|_{L^{p_i}}<\infty,$ then the function $\prod\limits^{\infty}_{i=1}f_i$ is well defined and $\|\prod\limits^{\infty}_{i=1}f_i\|_{L^1}\leq \prod\limits^{\infty}_{i=1}\|f_i\|_{L^{p_i}}.$
\end{lemma}

\noindent{\bf Proof of Lemma \ref{C-P-lem-G}} We split the proof into three cases.

Firstly, we assume that $\prod\limits^{\infty}_{i=1}\|f_i\|_{L^{p_i}}=0$ and there exists an $i_0\in N$ such that $\|f_{i_0}\|_{L^{p_i}}=0.$ It is clear that the function $\prod\limits^{\infty}_{i=1}f_i$ is well defined and $\|\prod\limits^{\infty}_{i=1}f_i\|_{L^1}\leq \prod\limits^{\infty}_{i=1}\|f_i\|_{L^{p_i}}.$

Secondly, we assume that $\prod\limits^{\infty}_{i=1}\|f_i\|_{L^{p_i}}=0$ and $\|f_i\|_{L^{p_i}}>0,~\forall i\in N.$
Let $\hat{f_i}=\frac{f_i}{\|f_i\|_{L^{p_i}}},~i\in N.$ Then $\|\hat{f_i}\|_{L^{p_i}}=1,~ i\in N.$
It follows from Lemma \ref{C-P-lem-F} that $\prod\limits^{\infty}_{i=1}\hat{f_i}$ is well defined. Combining this with $f_i=\|f_i\|_{L^{p_i}}\cdot\hat{f_i},~i\in N,$ we obtain that $\prod\limits^{\infty}_{i=1}f_i$ is well defined and $$\prod\limits^{\infty}_{i=1}f_i=\prod\limits^{\infty}_{i=1}\|f_i\|_{L^{p_i}}\prod\limits^{\infty}_{i=1}\hat{f_i}=0.$$
Thus, $\|\prod\limits^{\infty}_{i=1}f_i\|_{L^1}\leq \prod\limits^{\infty}_{i=1}\|f_i\|_{L^{p_i}}.$

Finally, we suppose that $0<\prod\limits^{\infty}_{i=1}\|f_i\|_{L^{p_i}}<\infty.$ Let $\hat{f_i}=\frac{f_i}{\|f_i\|_{L^{p_i}}},~i\in N.$ Then $\|\hat{f_i}\|_{L^{p_i}}=1,~ i\in N.$
It follows from Lemma \ref{C-P-lem-F} that $\prod\limits^{\infty}_{i=1}\hat{f_i}$ is well defined and $$\|\prod\limits^{\infty}_{i=1}\hat{f_i}\|_{L^1}\leq 1.$$
Thus the function $\prod\limits^{\infty}_{i=1}f_i$ is also well defined and $\|\prod\limits^{\infty}_{i=1}f_i\|_{L^1}\leq \prod\limits^{\infty}_{i=1}\|f_i\|_{L^{p_i}}.\blacksquare$

\begin{theorem}\label{C-P-prop-A}
Let $0<p_i<\infty, i\in N$ and $\sum\limits_{i=1}^{\infty}\frac{1}{p_i}=\frac{1}{p}.$
If $\prod\limits^{\infty}_{i=1}\|f_i\|_{L^{p_i}}<\infty,$ then the function $\prod\limits^{\infty}_{i=1}f_i$ is well defined and $\|\prod\limits^{\infty}_{i=1}f_i\|_{L^p}\leq \prod\limits^{\infty}_{i=1}\|f_i\|_{L^{p_i}}.$
\end{theorem}

\noindent{\bf Proof of Theorem \ref{C-P-prop-A}} It is clear that
Theorem \ref{C-P-prop-A} follows from Lemma \ref{C-P-lem-G}.$\blacksquare$

\begin{remark}Karakostas \cite{Karakostas} got the following result which was discussed on the $\sigma-$finite measure space. Let $1\leq p_i\leq\infty, i\in N$ and $\sum\limits_{i=1}^{\infty}\frac{1}{p_i}=\frac{1}{p}.$
If $0<\prod\limits^{\infty}_{i=1}\|f_i\|_{L^{p_i}}\leq\infty$ and the function $\prod\limits^{\infty}_{i=1}f_i$ is well defined, then $\|\prod\limits^{\infty}_{i=1}f_i\|_{L^p}\leq \prod\limits^{\infty}_{i=1}\|f_i\|_{L^{p_i}}.$
\end{remark}
\begin{theorem}\label{C-P-series-cor-E}
Let $1<p_i<\infty, i\in N$ and $\sum\limits_{i=1}^{\infty}\frac{1}{p_i}=\frac{1}{p}.$
Then
$$\prod\limits^{\infty}_{i=1}p'_i<\infty,$$
where $\frac{1}{p_i}+\frac{1}{p'_i}=1,~i\in N.$
\end{theorem}
\noindent{\bf Proof of Theorem \ref{C-P-series-cor-E} } It suffices to prove
$\sum\limits_{i=1}^{\infty}\ln p'_i<\infty.$ Because of $p'_i=(1-\frac{1}{p_i})^{-1},$ we should prove
$\sum\limits_{i=1}^{\infty}\ln (1-\frac{1}{p_i})^{-1}<\infty.$
Since $\lim\limits_{i\rightarrow\infty}\frac{\ln (1-\frac{1}{p_i})^{-1}}{\frac{1}{p_i}}=1$
and $\sum\limits_{i=1}^{\infty}\frac{1}{p_i}=\frac{1}{p},$
we have $\sum\limits_{i=1}^{\infty}\ln (1-\frac{1}{p_i})^{-1}<\infty$ by the Limit Comparison Test.$\blacksquare$

Let $\Omega=R^n$ and let $v_n$ be the volume of the unit ball in $R^n.$ If $f\in L^1,$ it follows from \cite[Theorem $2.1.6$]{L. Grafakos} and \cite[Exercise $2.1.3$]{L. Grafakos} that $\sup_{\lambda>0}\lambda|\{Mf>\lambda\}|\leq\xi_n\|f\|_{L^1},$ where $\xi_n=3^n\big(v_n(n/2)^{n/2}\big)^{-1}.$ It is a trivial fact that
$M$ maps $L^{\infty}\rightarrow L^{\infty}$ with constant $1.$ Using \cite[Exercise $1.3.3$]{L. Grafakos}, we obtain the following estimate
$$\|Mf\|_{L^p}\leq p'\xi_n^\frac{1}{p}\|f\|_{L^p}$$
for all $f\in L^p,~1<p<\infty.$ Then we have the following Lemma \ref{HL-constant}.

\begin{lemma}\label{HL-constant}Let $1<p_i<\infty,~i\in N.$ If $\prod\limits^{\infty}_{i=1}\|f_i\|_{L^{p_i}}<\infty$ and $\sum\limits_{i=1}^{\infty}\frac{1}{p_i}=\frac{1}{p},$
then $$ \|\mathfrak{M}(\overrightarrow f)\|_{L^p}\leq\|\prod\limits^{\infty}_{i=1}Mf_i\|_{L^p}
       \leq \xi_n^{\sum\limits_{i=1}^{\infty}\frac{1}{p_i}}\big(\prod\limits^{\infty}_{i=1}p'_i\big)\prod\limits_{i=1}^{\infty}\|f_i\|_{L^{p_i}}
       = \xi_n^{\frac{1}{p}} \big(\prod\limits^{\infty}_{i=1}p'_i\big)\prod\limits_{i=1}^{\infty}\|f_i\|_{L^{p_i}}<\infty.$$
\end{lemma}

\subsection{Generalized H\"{o}lder's Inequality for Conditional Expectation}\label{section B 2}

In the subsection,
we suppose that $(\Omega,\mathcal{F},\mu)$ is a complete probability space and
$\{f_i\}$ is a sequence of nonnegative measurable functions on $(\Omega,\mathcal{F},\mu).$

\begin{prop}\label{C-P-prop-gen-h}
Let $1<p_i<\infty, i\in N$ and $\sum\limits_{i=1}^{\infty}\frac{1}{p_i}=1.$
Suppose that $\mathcal{F}'$ be a sub$\hbox{-}\sigma\hbox{-}$field of $\mathcal{F}.$
If $\prod\limits^{\infty}_{i=1}\|f_i\|_{L^{p_i}}<\infty,$
then
$$E_{\mathcal{F}'}\big(\prod\limits^{\infty}_{i=1}f_i\big)
\leq \prod\limits^{\infty}_{i=1}E_{\mathcal{F}'}(f_i^{p_i})^{\frac{1}{p_i}}<\infty.$$
\end{prop}

\noindent{\bf Proof of Proposition \ref{C-P-prop-gen-h} }
Because of $\prod\limits^{\infty}_{i=1}\|f_i\|_{L^{p_i}}<\infty,$ it follows from Lemma \ref{C-P-lem-G} that
the function $\prod\limits^{\infty}_{i=1}f_i$ is well defined and $\|\prod\limits^{\infty}_{i=1}f_i\|_{L^1}\leq \prod\limits^{\infty}_{i=1}\|f_i\|_{L^{p_i}}.$ Since $\|f_i\|_{L^{p_i}}=\|f_i^{p_i}\|^{\frac{1}{p_i}}_{L^1}
=\|E_{\mathcal{F}'}(f_i^{p_i})\|^{\frac{1}{p_i}}_{L^1}=\|E_{\mathcal{F}'}(f_i^{p_i})^{\frac{1}{p_i}}\|_{L^{p_i}},~\forall i\in N,$ we have that
$\prod\limits^{\infty}_{i=1}E_{\mathcal{F}'}(f_i^{p_i})^{\frac{1}{p_i}}$ is well defined and
$$\|\prod\limits^{\infty}_{i=1}E_{\mathcal{F}'}(f_i^{p_i})^{\frac{1}{p_i}}\|_{L^1}
\leq\prod\limits^{\infty}_{i=1}\|f_i\|_{L^{p_i}}<\infty.$$
Moreover, $\prod\limits^{\infty}_{i=1}E_{\mathcal{F}'}(f_i^{p_i})^{\frac{1}{p_i}}<\infty.$ So we will focus on proving $E_{\mathcal{F}'}\big(\prod\limits^{\infty}_{i=1}f_i\big)
\leq \prod\limits^{\infty}_{i=1}E_{\mathcal{F}'}(f_i^{p_i})^{\frac{1}{p_i}}.$
For $k\in N,$ we define $q_k=\frac{1}{\sum\limits_{i=1}^{k}\frac{1}{p_i}},$ then $\sum\limits_{i=1}^{k}\frac{1}{p_i}=\frac{1}{q_k}$. Applying Fatou's Lemma and H\"{o}lder's inequality for conditional expectation, we have
\be E_{\mathcal{F}'}(\prod\limits^{\infty}_{i=1}f_i)
                &\leq&\liminf\limits_{k\rightarrow\infty}E_{\mathcal{F}'}(\prod\limits^{k}_{i=1}f_i)\\
                &\leq&\liminf\limits_{k\rightarrow\infty}E_{\mathcal{F}'}(\prod\limits^{k}_{i=1}f^{q_k}_i)^{\frac{1}{q_k}}\\
                &\leq&\liminf\limits_{k\rightarrow\infty}\prod\limits^{k}_{i=1}E_{\mathcal{F}'}(f^{p_i}_i)^{\frac{1}{p_i}}\\
                &=&\prod\limits^{\infty}_{i=1}E_{\mathcal{F}'}(f_i^{p_i})^{\frac{1}{p_i}}.\blacksquare\ee

\begin{theorem}\label{C-P-cor-gen-h}
Let $0<p_i<\infty, i\in N$ and $\sum\limits_{i=1}^{\infty}\frac{1}{p_i}=\frac{1}{p}.$
Suppose that $\mathcal{F}'$ be a sub$\hbox{-}\sigma\hbox{-}$field of $\mathcal{F}.$
If $\prod\limits^{\infty}_{i=1}\|f_i\|_{L^{p_i}}<\infty,$
then $$E_{\mathcal{F}'}\big(\prod\limits^{\infty}_{i=1}f_i^p\big)^{\frac{1}{p}}
\leq \prod\limits^{\infty}_{i=1}E_{\mathcal{F}'}(f_i^{p_i})^{\frac{1}{p_i}}.$$
\end{theorem}

\noindent{\bf Proof of Theorem \ref{C-P-cor-gen-h}} It is clear that
Theorem \ref{C-P-cor-gen-h} follows from Proposition \ref{C-P-prop-gen-h}.$\blacksquare$

\begin{prop} \label{C-P-lem-def} Let $1<p_i<\infty, i\in N$ and $\sum\limits_{i=1}^{\infty}\frac{1}{p_i}=\frac{1}{p}.$
If $\prod\limits_{i=1}^{\infty}\|f_i\|_{L^{p_i}}<\infty,$ then $\mathfrak{M}(\overrightarrow{f})$ is well defined.
\end{prop}

\noindent{\bf Proof of Proposition \ref{C-P-lem-def}} Let $q>1.$ It is well known that conditional
expectation $E_n(\cdot)$ on $L^q(\Omega,\mathcal{F},\mu)$ is a contraction,
and maps $L^q(\Omega,\mathcal{F},\mu)$ onto $L^q(\Omega,\mathcal{F}_n,\mu).$
Combining this with Theorem \ref{C-P-prop-A} and Remark \ref{re-aa},
we have $\prod\limits^{\infty}_{i=1}E_n(f_i)$ is well defined.
Then $\mathfrak{M}(\overrightarrow{f})$ is well defined.$\blacksquare$

\section{Weighted Inequalities In Martingale Spaces}\label{section C}
There are several assumptions that will be used in this section.
For convenience, we state them at the beginning of this part. In addition, $C$ will denote a
constant not necessarily the same at each occurrence.

\textbf{ASSUMPTIONS} Let $\omega_i\in L^1$ and $1<p_i<\infty$, $i\in N,$ and let $\{f_i\}$ be a sequence of nonnegative measurable function on the probability space $(\Omega,\mathcal{F},\mu).$
Suppose that
$\frac{1}{p}=\sum\limits^{\infty}_{i=1}\frac{1}{p_i }$ and $\sigma_i=\omega_i^{-\frac{1}{p_i-1}}\in L^1,~i\in N.$
We always suppose that $\prod\limits_{i=1}^{\infty}\|\sigma_i\|_{L^{p_i}(\omega_i)}<\infty,$ $\prod\limits_{i=1}^{\infty}E_n(\omega_i^{1-p'_i})^{\frac{1}{p'_i}}<\infty,$
and $\prod\limits_{i=1}^{\infty}\|f_i\|_{L^{p_i}(\omega_i)}<\infty.$
Moreover, we assume that $\prod\limits_{i=1}^{\infty}\sigma_i^{\frac{1}{p_i}}>0.$

\textbf{NOTATIONS} We denote that $\overrightarrow{p}=(p_1,p_2,\cdot\cdot\cdot),$
$\overrightarrow{\omega}=(\omega_1, \omega_2,\cdot\cdot\cdot),$
$\overrightarrow{f}=(f_1,f_2,...).$ Moreover, we also denote $\overrightarrow{f\chi_Q}=(f_1\chi_Q,f_2\chi_Q,...)$
and $\overrightarrow{\sigma\chi_Q}
=(\sigma_1\chi_Q,\sigma_2\chi_Q,\cdot\cdot\cdot),$ where $Q$ is a measurable set.

\begin{remark} \label{C-P-lem-finite} It follows from the generalized H\"{o}lder's inequality for integral that $$\int_\Omega\prod\limits_{i=1}^{\infty}E_n(f_i^{p_i}\omega_i)^{\frac{p}{p_i}}d\mu
\leq\prod\limits_{i=1}^{\infty}\big(\int_\Omega E_n(f_i^{p_i}\omega_i)d\mu\big)^{\frac{p}{p_i}}
=\prod\limits_{i=1}^{\infty}\big(\int_\Omega f_i^{p_i}\omega_id\mu\big)^{\frac{p}{p_i}}<\infty.$$
Hence,
$\prod\limits_{i=1}^{\infty}E_n(f_i^{p_i}\omega_i)^{\frac{1}{p_i}}<\infty.$
By H\"{o}lder's inequality for conditional expectation and Remark \ref{re-aa}, we have
\be\prod\limits_{i=1}^{\infty}E_n(f_i)
   &\leq&\prod\limits_{i=1}^{\infty}E_n(f_i^{p_i}\omega_i)^{\frac{1}{p_i}}
         E_n(\omega_i^{-\frac{1}{p_i-1}})^{\frac{1}{p'_i}}\\
   &=&\prod\limits_{i=1}^{\infty}E_n(f_i^{p_i}\omega_i)^{\frac{1}{p_i}}
         \prod\limits_{i=1}^{\infty}E_n(\omega_i^{-\frac{1}{p_i-1}})^{\frac{1}{p'_i}}<\infty.\ee
Then $\mathfrak{M}(\overrightarrow{f})$ is well defined.
Let $f_i=\sigma_i,$
we also have $\prod\limits_{i=1}^{\infty}E_n(\sigma_i)<\infty$
and $\mathfrak{M}(\overrightarrow{\sigma})$ is well defined.\end{remark}

\begin{definition} \label{Rh}
We say that the weight vector $\overrightarrow{\omega}$
satisfies the reverse H\"{o}lder's condition $RH_{\overrightarrow{p}},$ if
there exists a positive constant $C$ such that
\be\prod\limits_{i=1}^{\infty}\big(\int_{\{\tau<\infty\}}\sigma_id\mu\big)^{\frac{p}{p_i}}
\leq C\int_{\{\tau<\infty\}}\prod\limits_{i=1}^{\infty}\sigma_i^{\frac{p}{p_i}}d\mu,~\forall \tau\in\mathcal{T}.\ee
\end{definition}

\begin{definition}Let $v$ be a weight.
We say that the weight vector $(v,\overrightarrow{\omega})$
satisfies the condition $A_{\overrightarrow{p}},$ if
there exists a positive constant $C$ such that
$$E_n(v)^{\frac{1}{p}}\prod\limits_{i=1}^{\infty}E_n(\omega_i^{1-p'_i})^{\frac{1}{p'_i}}
      \leq C,~\forall n\geq0,$$
where $\frac{1}{p_i}+\frac{1}{p'_i}=1,~i\in N.$
\end{definition}

\begin{definition}\label{Sp}Let $v$ be a weight.
We say that the weight vector $(v,\overrightarrow{\omega})$
satisfies the condition $S_{\overrightarrow{p}},$ if
there exists a positive constant $C$ such that
\be
\big(\int_{\{\tau<\infty\}}\mathfrak{M}(\overrightarrow{\sigma\chi_{\{\tau<\infty\}}}
)^{p}vd\mu\big)^{\frac{1}{p}}\leq
C\prod\limits_{i=1}^{\infty}|\{\tau<\infty\}|^{\frac{1}{p_i}}_{\sigma_i},
~\forall \tau\in \mathcal{T}.\ee
\end{definition}

\noindent{\bf Proof of Theorem \ref{thm_AP} } We shall follow the scheme:
$\eqref{C-L-P-W-2}\Leftrightarrow\eqref{C-L-P-W-1}\Leftrightarrow\eqref{C-L-P-W-3}.$

$(\ref{C-L-P-W-1})\Rightarrow(\ref{C-L-P-W-2}).$ Let $f_i\in L^{p_i}(\omega_i),~i\in N$ and let $\prod\limits_{i=1}^{\infty}\|f_i\|_{L^{p_i}(\omega_i)}<\infty.$ For
$\lambda>0,$ define $\tau=\inf\{n:\prod\limits_{i=1}^{\infty}E_n(f_i)>\lambda\}.$ It follows from $\eqref{Th_B_1}$ that
\be\lambda|\{\mathfrak{M}(\overrightarrow f)>\lambda\}|_v^{\frac{1}{p}}
   &=&(\int_{\{\tau<\infty\}}\lambda^p vd\mu )^{\frac{1}{p}}\\
   &\leq&\Big(\int_{\{\tau<\infty\}}
      \prod\limits_{i=1}^{\infty}E_{\tau}(f_i)^pvd\mu\Big)^{\frac{1}{p}}\\
   &\leq& C\prod\limits_{i=1}^{\infty}\|f_i\|_{L^{p_i}(\omega_i)}.
   \ee
Thus $\eqref{Th_B_2}$ is valid.

$(\ref{C-L-P-W-2})\Rightarrow(\ref{C-L-P-W-1}).$ Let $f_i\in L^{p_i}(\omega_i),~i\in N$ and let $\prod\limits_{i=1}^{\infty}\|f_i\|_{L^{p_i}(\omega_i)}<\infty.$
Fix $n\in N$ and $B\in\mathcal {F}_n.$ Let $$F_i=f_i\chi_B,~i \in N.$$
Then $E_n(F_i)=E_n(f_i)\chi_B.$
Moreover
$$
\prod\limits_{i=1}^{\infty}E_n(f_i)\chi_B\leq\mathfrak{M}(\overrightarrow F).
$$
Combining with $\eqref{Th_B_2},$ we have \be\lambda^p
        \int_{B\cap\{\prod\limits_{i=1}^{\infty}E_n(f_i)>\lambda\}}vd\mu
&\leq& \lambda^p\int_{\{\mathfrak{M}(\overrightarrow F)>\lambda\}}vd\mu \\
&\leq& C\prod\limits_{i=1}^{\infty}\|F_i\|^p_{L^{p_i}(\omega_i)}\\
&=& C\prod\limits_{i=1}^{\infty}\Big(\int_B f_i^{p_i}\omega_i
        d\mu\Big)^{\frac{p}{p_i}}.\ee
For $k\in Z$, let
\be
B_k = \{2^k<
    \prod\limits_{i=1}^{\infty}E_n(f_i)\leq2^{k+1}\}.
\ee
Note that $$\{2^k<\prod\limits_{i=1}^{\infty}E_n(f_i)\leq2^{k+1}\}
\subseteq \{2^k<\prod\limits_{i=1}^{\infty}E_n(f_i)\},$$
then
\be \int_\Omega
    (\prod\limits_{i=1}^{\infty}E_n(f_i))^pvd\mu
&=& \sum\limits_{k\in Z}\int_{B_k}
    (\prod\limits_{i=1}^{\infty}E_n(f_i))^pvd\mu\\
&\leq&C \sum\limits_{k\in Z}\int_{{B_k}\cap\{
    \prod\limits_{i=1}^{\infty}E_n(f_i)>2^k\}}2^{kp}vd\mu\\
&\leq&C \sum\limits_{k\in Z}\prod\limits_{i=1}^{\infty}\Big(\int_{B_k}f_i^{p_i}\omega_i
        d\mu\Big)^{\frac{p}{p_i}}\\
&\leq&C \prod\limits_{i=1}^{\infty}\Big(\sum\limits_{k\in Z}\int_{B_k}f_i^{p_i}\omega_i
        d\mu\Big)^{\frac{p}{p_i}}\\
&\leq&C\prod\limits_{i=1}^{\infty}\Big(\int_\Omega f_i^{p_i}\omega_i
        d\mu\Big)^{\frac{p}{p_i}},\ee
where we have used the generalized H\"{o}lder's inequality.
As for $\tau\in\mathcal {T},$
it is easy to see that
\be\int_{\{\tau<\infty\}}\prod\limits_{i=1}^{\infty}E_{\tau}(f_i)^pvd\mu
   &=&\sum\limits_{n\geq0}\int_{\{\tau=n\}}\prod\limits_{i=1}^{\infty}E_n(f_i)^pvd\mu\\
   &\leq&C\sum\limits_{n\geq0}\prod\limits_{i=1}^{\infty}\Big(\int_\Omega(f_i\chi_{\{\tau=n\}})^{p_i}\omega_i
        d\mu\Big)^{\frac{p}{p_i}}\\
   &\leq&C\prod\limits_{i=1}^{\infty}\Big(\sum\limits_{n\geq0}\int_\Omega (f_i\chi_{\{\tau=n\}})^{p_i}\omega_i
        d\mu\Big)^{\frac{p}{p_i}}\\
   &\leq&C\prod\limits_{i=1}^{\infty}\Big(\int_\Omega f_i^{p_i}\omega_i
        d\mu\Big)^{\frac{p}{p_i}}.\ee
Therefore, $$\Big(\int_{\{\tau<\infty\}}
      \prod\limits_{i=1}^{\infty}E_{\tau}(f_i)^pvd\mu\Big)^{\frac{1}{p}}
      \leq C\prod\limits_{i=1}^{\infty}\|f_i\|_{L^{p_i}(\omega_i)}.$$

$\eqref{C-L-P-W-3}\Rightarrow\eqref{C-L-P-W-1}.$ Let $f_i\in L^{p_i}(\omega_i),~i\in N$ and let $\prod\limits_{i=1}^{\infty}\|f_i\|_{L^{p_i}(\omega_i)}<\infty.$
Applying H\"{o}lder's inequality for conditional expectation, we get
$$E_n(f_i)\leq E_n(f_i^{p_i}\omega_i)^{\frac{1}{p_i}}E_n(\omega_i^{-\frac{1}{p_i-1}})^{\frac{1}{p'_i}}.$$
Furthermore,
\be\prod\limits_{i=1}^{\infty}E_n(f_i)^p
   &\leq&\prod\limits_{i=1}^{\infty}E_n(f_i^{p_i}\omega_i)^{\frac{p}{p_i}}
   E_n(\omega_i^{-\frac{1}{p_i-1}})^{\frac{p}{p'_i}}\\
&=&\prod\limits_{i=1}^{\infty}E_n^v(f^{p_i}\omega_iv^{-1})^{\frac{p}{p_i}}
   E_n(v)E_n(\omega_i^{-\frac{1}{p_i-1}})^{\frac{p}{p'_i}},\ee
where $E_n^v(\cdot)$ is the conditional expectation relative to the probability measure $\frac{v}{|\Omega|_v}d\mu.$
Because of $(\ref{Th_B_3}),$ we get $$\prod\limits_{i=1}^{\infty}E_n(f_i)^p
   \leq C\prod\limits_{i=1}^{\infty}E_n^v(f^{p_i}\omega_iv^{-1})^{\frac{p}{p_i}}.$$
From this, using the generalized H\"{o}lder's inequality, we have
\be\|\prod\limits_{i=1}^{\infty}E_n(f_i)\|_{L^p(v)}
   &\leq& C\|\prod\limits_{i=1}^{\infty}E_n^v(f^{p_i}\omega_iv^{-1})^{\frac{1}{p_i}}\|_{L^p(v)}\\
   &\leq& C\prod\limits_{i=1}^{\infty}\|E_n^v(f^{p_i}\omega_iv^{-1})^{\frac{1}{p_i}}\|_{L^{p_i}(v)}\\
   &=& C\prod\limits_{i=1}^{\infty}\|E_n^v(f^{p_i}\omega_iv^{-1})\|^{\frac{1}{p_i}}_{L^{1}(v)}\\
   &=& C\prod\limits_{i=1}^{\infty}\|f^{p_i}\omega_i\|^{\frac{1}{p_i}}_{L^{1}}\\
   &=& C\prod\limits_{i=1}^{\infty}\|f^{p_i}\|_{L^{p_i}(\omega_i)}.\ee

$(\ref{C-L-P-W-1})\Rightarrow(\ref{C-L-P-W-3}).$ For any $n\geq 0,$ $i\in N$ and $B\in \mathcal{F}_n,$ set
$f_i=\omega_i^{-\frac{1}{p_i-1}}\chi_B.$
Then $$\Big(\int_B\prod\limits_{i=1}^{\infty}E_n(\omega_i^{-\frac{1}{p_i-1}})^pvd\mu\Big)^{\frac{1}{p}}
      \leq C\prod\limits_{i=1}^{\infty}\Big(\int_\Omega \omega_i^{-\frac{1}{p_i-1}}\chi_Bd\mu\Big)^{\frac{1}{p_i}}.$$
Furthermore,
$$\int_B\prod\limits_{i=1}^{\infty}
E_n(\omega_i^{-\frac{1}{p_i-1}})^pE_n(v)d\mu
      \leq C\prod\limits_{i=1}^{\infty}\Big(\int_B
      \sigma_id\mu\Big)^{\frac{p}{p_i}}.$$
Note that $\overrightarrow{\omega}\in RH_{\overrightarrow{p}},$ we have $$\int_B\prod\limits_{i=1}^{\infty}
E_n(\omega_i^{-\frac{1}{p_i-1}})^pE_n(v)d\mu
      \leq C\int_B
      \prod\limits_{i=1}^{\infty}\sigma_i^{\frac{p}{p_i}}d\mu.$$
It follows from the generalized H\"{o}lder's inequality for conditional expectation that $$\int_B\prod\limits_{i=1}^{\infty}
E_n(\omega_i^{-\frac{1}{p_i-1}})^pE_n(v)d\mu
      \leq C\int_B
      \prod\limits_{i=1}^{\infty}E_n(\sigma_i)^{\frac{p}{p_i}}d\mu.$$
Thus, there exists a constant $C$ such that$$\Big(\prod\limits_{i=1}^{\infty}
E_n(\omega_i^{-\frac{1}{p_i-1}})^pE_n(v)\Big)^{\frac{1}{p}}
      \leq C\prod\limits_{i=1}^{\infty}
      E_n(\omega_i^{-\frac{1}{p_i-1}})^{\frac{1}{p_i}}.$$
Then
\be E_n(v)^{\frac{1}{p}}\prod\limits_{i=1}^{\infty}E_n(\omega_1^{1-p'_i})^{\frac{1}{p'_i}}\leq C.\blacksquare\ee

\noindent{\bf Proof of Theorem \ref{thm_Sp}} It is clear that $\eqref{C_L_P_thm_Sp_1}\Leftrightarrow
\eqref{C_L_P_thm_Sp_2}\Rightarrow\eqref{C_L_P_thm_Sp_3},$
so we omit them.  To prove $\eqref{C_L_P_thm_Sp_3}\Rightarrow\eqref{C_L_P_thm_Sp_2},$
we proceed in the following way. Let $g_i\in L^{p_i}(\sigma_i),~i\in N$ and let $\prod\limits_{i=1}^{\infty}\|g_i\|_{L^{p_i}(\sigma_i)}<\infty.$
For all $k\in Z$, define
stopping times $$\tau_k=\inf\{n:\prod\limits_{i=1}^{\infty}E_n(g_i\sigma_i)>2^k\}.$$
Set
\be
A_{k,j}=\{\tau_k<\infty\}\cap\{2^j<\prod\limits_{i=1}^{\infty}E_{\mathcal{F}_{\tau_k}}(\sigma_i)\leq2^{j+1}\};
\ee
\be
B_{k,j}=\{\tau_k<\infty,\tau_{k+1}=\infty\}\cap\{2^j<\prod\limits_{i=1}^{\infty}
E_{\mathcal{F}_{\tau_k}}(\sigma_i)\leq2^{j+1}\},~ j\in Z.
\ee
Then $A_{k,j}\in \mathcal{F}_{\tau_k}, B_{k,j}\subseteq A_{k,j}$ and
\be
E_{\mathcal{F}_{\tau_k}}(g_i\sigma_i)=E^{\sigma_i}_{\mathcal{F}_{\tau_k}}(g_i)E_{\mathcal{F}_{\tau_k}}(\sigma_i).\ee
Moreover, $\{B_{k,j}\}_{k,j}$ is a family of disjoint sets and
\be
\{2^k<\mathfrak{M}(\overrightarrow {g\sigma})\leq2^{k+1}\}=\{\tau_k<\infty,\tau_{k+1}=\infty\}=\bigcup\limits_{j\in
Z} B_{k,j}, k\in Z.
\ee

On each $A_{k,j},$ we have \be 2^{kp}&\leq&
                 \mathop{\hbox{ess inf}}\limits_{A_{k,j}}
                 \prod\limits_{i=1}^{\infty}E_{\mathcal{F}_{\tau_k}}(g_i\sigma_i)^p\\
              &\leq&\mathop{\hbox{ess inf}}\limits_{A_{k,j}}
                 \prod\limits_{i=1}^{\infty}E^{\sigma_i}_{\mathcal{F}_{\tau_k}}(g_i)^p
                 \mathop{\hbox{ess sup}}\limits_{A_{k,j}}
                 \prod\limits_{i=1}^{\infty}E_{\mathcal{F}_{\tau_k}}(\sigma_i)^p\\
              &\leq&2^p\mathop{\hbox{ess inf}}\limits_{A_{k,j}}
                 \prod\limits_{i=1}^{\infty}E^{\sigma_i}_{\mathcal{F}_{\tau_k}}(g_i)^p
                 |B_{k,j}|_v^{-1}\int_{B_{k,j}}
                 \prod\limits_{i=1}^{\infty}E_{\mathcal{F}_{\tau_k}}(\sigma_i)^pvd\mu.\ee
To estimate $\int_\Omega\mathfrak{M}(\overrightarrow {g\sigma})^p v d\mu, $ firstly we have \be
\int_\Omega\mathfrak{M}(\overrightarrow {g\sigma})^p v d\mu
   &=&\sum\limits_{k\in Z}\int_{\{2^k<\mathfrak{M}(\overrightarrow {g\sigma})\leq2^{k+1}\}}
      \mathfrak{M}(\overrightarrow {g\sigma})^p v d\mu\\
   &\leq&2^p\sum\limits_{k\in Z}\int_{\{2^k<\mathfrak{M}(\overrightarrow {g\sigma})\leq2^{k+1}\}}2^{kp}vd\mu\\
   &=&2^p\sum\limits_{k\in Z,j\in Z}2^{kp}\int_{B_{k,j}} v d\mu\\
   &\leq&4^p\sum\limits_{k\in Z,j\in Z}\mathop{\hbox{ess inf}}\limits_{A_{k,j}}
                 \prod\limits_{i=1}^{\infty}E^{\sigma_i}_{\mathcal{F}_{\tau_k}}(g_i)^p
                 \int_{B_{k,j}}\prod\limits_{i=1}^{\infty}E_{\mathcal{F}_{\tau_k}}(\sigma_i)^pvd\mu.
\ee
It is clear that $\vartheta$ is a measure on $X=Z^2$ with
\be
\vartheta(k,j)=\int_{B_{k,j}}\prod\limits_{i=1}^{\infty}E_{\mathcal{F}_{\tau_k}}(\sigma_i)^pvd\mu.
\ee
For the above $\{g_i\},$ define \be
T_{\overrightarrow g}(k,j)=\mathop{\hbox{ess inf}}\limits_{A_{k,j}}
                 \prod\limits_{i=1}^{\infty}E^{\sigma_i}_{\mathcal{F}_{\tau_k}}(g_i)^p\ee
and denote
\be E_\lambda=\Big\{(k,j):
\mathop{\hbox{ess inf}}\limits_{A_{k,j}}
                 \prod\limits_{i=1}^{\infty}E^{\sigma_i}_{\mathcal{F}_{\tau_k}}(g_i)^p>\lambda\Big\}
\hbox{ and }G_\lambda=\bigcup\limits_{(k,j)\in E_\lambda}A_{k,j} \ee
for each $\lambda>0. $ Then we have \be |\{T_{\overrightarrow g}(k,j)>\lambda\}|_\vartheta
&=&\sum\limits_{(k,j)
        \in E_\lambda}\int_{B_{k,j}}\prod\limits_{i=1}^{\infty}E_{\mathcal{F}_{\tau_k}}(\sigma_i)^pvd\mu\\
&=&\sum\limits_{(k,j)
        \in E_\lambda}\int_{B_{k,j}}\prod\limits_{i=1}^{\infty}E_{\mathcal{F}_{\tau_k}}(\sigma_i\chi_{G_\lambda})^pvd\mu\\
&\leq& \int_{G_\lambda}\mathfrak{M}(\overrightarrow{\sigma\chi_{G_\lambda}})^p v d\mu. \ee
Let $\tau=\inf\Big\{n:~
          \prod\limits_{i=1}^{\infty}E^{\sigma_i}_n(g_i)^p>\lambda\Big\},$ we have
$G_\lambda\subseteq\Big\{\mathfrak{M}^{\overrightarrow\sigma}(\overrightarrow g)^p>\lambda\Big\}=\{\tau<\infty\}.$
It follows from $S_{\overrightarrow{p}}$ and $RH_{\overrightarrow{p}}$ that
 \be |\{T_{\overrightarrow g}(k,j)>\lambda\}|_\vartheta
&\leq&\int_{\{\tau<\infty\}}\mathfrak{M}(\overrightarrow{\sigma\chi_{\{\tau<\infty\}}})^p v d\mu.\\
&\leq& C\prod\limits_{i=1}^{\infty}|\{\tau<\infty\}|^{\frac{p}{p_i}}_{\sigma_i}\\
&\leq& C\int_{\{\tau<\infty\}}\prod\limits_{i=1}^{\infty}\sigma_i^{\frac{p}{p_i}}d\mu.
 \ee
Therefore,
\be~~~~~~\int_\Omega\mathfrak{M}(\overrightarrow {g\sigma})^p v d\mu
&\leq&4^p\int_XT_{\overrightarrow g}d\vartheta=4^p\int_0^\infty|\{T_{\overrightarrow g}>\lambda\}|_\vartheta
              d\lambda\\
&\leq&C\int_0^\infty\int_{\{\tau<\infty\}}\prod\limits_{i=1}^{\infty}\sigma_i^{\frac{p}{p_i}}d\mu d\lambda \\
&=&C\int_0^\infty\int_{\{\mathfrak{M}^{\overrightarrow\sigma}(\overrightarrow g)^p>\lambda\}}
              \prod\limits_{i=1}^{\infty}\sigma_i^{\frac{p}{p_i}}d\mu d\lambda \\
&=&C\int_\Omega\mathfrak{M}^{\overrightarrow\sigma}(\overrightarrow g)^p\prod\limits_{i=1}^{\infty}\sigma_i^{\frac{p}{p_i}}d\mu\\
&\leq&C\int_\Omega \prod\limits_{i=1}^{\infty} M^{\sigma_i}(g_i)^p\sigma_i^{\frac{p}{p_i}}d\mu\\
&\leq&C\prod\limits_{i=1}^{\infty}\big(\int_\Omega M^{\sigma_i}(g)^{p_i}\sigma_id\mu\big)^{\frac{p}{p_i}}\\
&\leq&C(\prod\limits_{i=1}^{\infty}p_i')^p\prod\limits_{i=1}^{\infty}\|g_i\|_{L^{p_i}(\sigma_i)}^{p},
\ee where we use H\"{o}lder's inequality and Doob's inequality.
Then $(\ref{Th_A_1})$ is valid, because of $\prod\limits_{i=1}^{\infty}p_i'<\infty.\blacksquare$

\thanks{\textbf{Acknowledgement} This paper was partially completed while W. Chen was at the Institute of Mathematics of the
University of Seville (I.M.U.S.), Spain. He would like to express his gratitude for the hospitality received there.
The authors thank the anonymous referee for his/her careful reading of the manuscript and useful corrections.
They also thank Peide Liu, Youliang Hou and Maofa Wang for many valuable comments on this paper.}

%
%

\end{document}